\documentclass[psamsfont,a4paper,12pt]{amsart}
\usepackage[english]{babel}
\numberwithin{equation}{section}

\usepackage{setspace}
\usepackage{enumerate}
\usepackage{afterpage}
\usepackage[usenames,dvipsnames]{color}
\usepackage{graphicx}

\usepackage{amssymb,latexsym,amsmath}
\usepackage{amssymb}
\usepackage{mathrsfs}
\usepackage{amsfonts}
\usepackage{latexsym}
\usepackage[latin1]{inputenc}
\usepackage{pstricks,multido,pst-plot}
\usepackage[none]{hyphenat}

\usepackage{amsthm}
\usepackage{stackrel}                         

\usepackage{tabularx}
\usepackage{longtable}

\usepackage{verbatim}
\usepackage{blindtext}
\usepackage{multicol,fancyheadings}
\usepackage{float}
\usepackage[vcentermath,enableskew]{youngtab}
\usepackage{subfig}

\allowdisplaybreaks

\newtheorem{theorem}{Theorem}[section]			%%with section numbering

\newtheorem{lemma}[theorem]{Lemma}

\theoremstyle{definition}

\theoremstyle{remark}

\usepackage{tikz}
\usetikzlibrary{arrows}

\usepackage{paralist}
\usepackage{cite}

\title{Non-standard quaternary representations and the Fibonacci numbers}

\author[Anders]{Katie Anders}
\address{Department of Mathematics, University of Texas at Tyler, Tyler, TX 75799}
\email{kanders@uttyler.edu}

\author[Dawsey]{Madeline L. Dawsey}
\address{Department of Mathematics, University of Texas at Tyler, Tyler, TX 75799}
\email{mdawsey@uttyler.edu}

\author[Gupta]{Rajat Gupta}
\address{Department of Mathematics \& Statistics, University of Maine, Orono, ME 04469}
\email{rajat.gupta@maine.edu}

\author[Lebowitz-Lockard]{Noah Lebowitz-Lockard}
\address{Department of Mathematics, University of Texas at Tyler, Tyler, TX 75799}
\email{nlebowi@gmail.com}

\author[Vandehey]{Joseph Vandehey}
\address{Department of Mathematics, University of Texas at Tyler, Tyler, TX 75799}
\email{jvandehey@uttyler.edu}


\begin{document}

%\keywords{}

\begin{abstract}
Let $f_4(n)$ be the number of hyperquaternary representations of $n$ and $b_4(n)$ be the number of balanced quaternary representations of $n$. We show that there is no integer $k$ such that $f_4(n+k)=b_4(n)$ for all $n\ge -k$, in contrast to the binary case. Nevertheless, there do exist integers $k$ such that $f_4(n+k)=b_4(n)$ for arbitrarily large intervals of $n$. We generalize these results to any even base $d$.  We also study the rate of growth of $b_4(n)$ and show that maximal values of this function correspond to certain Fibonacci numbers.
\end{abstract}

\maketitle

\section{Introduction}\label{introduction}

For integers $d\ge1 $ and $n\ge0$, a standard $d$-ary representation of $n$ is an expression for $n$ of the form
\begin{equation*}
n=\sum_{i=0}^\infty\epsilon_id^i,\quad\epsilon_i\in\{0,1,2,\dots,d-1\}.
\end{equation*}
Since almost all digits $\epsilon_i$ will be $0$, we will condense notation by writing
\begin{equation}\label{eq: digits representation}
[\epsilon_k \; \epsilon_{k-1} \; \cdots \; \epsilon_0]_d = \sum_{i=0}^k \epsilon_id^i.
\end{equation}
To further simplify lengthy representations, we will use two additional notations. If a single digit is repeated several times, we will write a brace underneath with a number to indicate the total number of digits present. For example, $[\underbrace{1 \; 1\; \dots \; 1}_j]_d$ indicates that there are $j$ $1$'s in this representation. If there are repetitions of a combination of digits, we will denote this with a power indicating the number of repetitions. As an example, $[(1\; 0)^3]_d=[1\;0 \; 1\; 0 \;1 \; 0]_d$.

It is known that standard $d$-ary representations are unique for each $n\in\mathbb{N}_{\ge 0}$.  If the digit set $\{0,1,2,\dots,d-1\}$ is modified, such modified $d$-ary representations may not be unique but may exhibit interesting number theoretic properties and are connected with the study of restricted partition functions. For more information, see, for example, \cite{anders2013congruence,colson1726short,hare2015base,protasov2000asymptotic,reznick1990some,tangjai2020non}.

In this paper, we consider two types of non-standard $d$-ary representations. The first, the hyper-$d$-ary representation, uses the digit set $\{0,1,2,\dots, d\}$, while the second, the balanced $d$-ary representation, uses the digit set $\mathbb{Z} \cap [-d/2, d/2]$. If $d$ is odd, there is a unique  balanced $d$-ary representation for each integer; whereas for $d$ even, some integers have multiple balanced $d$-ary representations.  The fact that balanced ternary ($d=3$) gives unique representations was known to Euler and even earlier mathematicians (see \cite{andrews2007euler}); however, the uniqueness of representations for odd $d$ and non-uniqueness of representations for even $d$ are both simple to prove. For example, in balanced binary, there are infinitely many representations of any given integer, such as $1=[1]_2=[1\; -\!1]_2 = [1\; -\!1\; -\!1]_2=\cdots$.

We let $f_d(n)$ denote the number of hyper-$d$-ary representations of $n$ and $b_d(n)$ denote the number of balanced $d$-ary representations of $n$. Reznick \cite{reznick1990some} showed that $f_2(n)=s(n+1)$, where $s(n)$ is the Stern sequence, given by $s(0)=0$, $s(1)=1$, and 
\begin{equation}
s(2n)= s(n), \qquad s(2n+1)=s(n)+s(n+1)\label{eq:stern recurrence}
\end{equation} for $n\geq2$. In a previous paper \cite{anders2023non}, several of the authors of this paper studied $\overline{b}_2(n)$, which counts the number of \emph{short} balanced binary representations of $n$.  These are balanced binary representations not of the form $[1\; -\!1 \; \cdots]_2$ or $[-1 \; 1 \; \cdots]_2$.\footnote{The authors opted to study $\overline{b}_2(n)$ instead of $b_2(n)$ as $b_2(n)$ is trivially infinite for all positive $n$, by the remark at the end of the previous paragraph.} The authors showed that $\overline{b}_2(n)=s(n)$ via several means, including a direct proof, a bijective proof, and a generating function proof.  Combining this result with Reznick's, we have that $\overline{b}_2(n)=f_2(n-1)$ for $n\ge 1$. This leads to the natural question: for other even bases $d$, does there exist an integer $k$ such that $b_d(n)=f_d(n+k)$? 

At present, we restrict our attention to the case $d=4$, and then we will return to the general case in Section \ref{sec:higher bases}. We now consider two types of non-standard $4$-ary, or \emph{quaternary}, representations, which are special cases of the representations described earlier.  A \emph{hyperquaternary} representation of $n$ is an expression for $n$ of the form
\begin{equation}\label{hyperquadrary_definition}
n=\sum_{i=0}^\infty\epsilon_i4^i,\quad\epsilon_i\in\{0,1,2,3,4\}.
\end{equation}
Note that $f_4(n)$ counts the number of hyperquaternary representations of an integer $n$, and observe that $f_4(n)=0$ when $n<0$.  For example, $f_4(67)=3$ because $67$ can be written in hyperquaternary as $[1 \; 0 \; 0 \; 3]_4$, as $[4\; 0 \; 3]_4$, and as $[3\; 4\; 3]_4$.

A \emph{balanced quaternary} representation of $n$ is an expression for $n$ of the form
\begin{equation}\label{balanced_quaternary_definition}
n=\sum_{i=0}^\infty\epsilon_i4^i,\quad\epsilon_i\in\{-2,-1,0,1,2\}.
\end{equation}
 Note that for any integer $n$, $b_4(n)$ counts the number of balanced quaternary representations of $n$.  For example, $b_4(25)=3$ because $25$ can be written in balanced quaternary as $[1 \; 2 \; 1]_4$, as $[2 \; -\!2 \; 1]_4$, and as $[1 \; -\!2 \; -\!2 \; 1]_4$.

Our first result answers our question above in the negative. 

\begin{theorem}\label{thm:no good k}
    There does not exist an integer $k$ such that $f_4(n+k)=b_4(n)$ for all $n\ge -k$.
\end{theorem}

Nevertheless, we can find a $k$ such that $f_4(n+k)=b_4(n)$ for a long interval of $n$'s, as we show in the following result.

\begin{theorem}\label{thm:equality on interval}
    For each $j\in \mathbb{N}$, we have that
    \[
     f_4\big([\underbrace{2 \; 2\; \cdots \; 2}_{j}]_4+n\big)=b_4(n), \qquad n \in \left[ -[\underbrace{1\; 1\; \cdots \; 1}_{j}]_4 , [\underbrace{1\; 1\; \cdots \; 1}_{j+1}]_4 \right].
    \]
    Moreover, $f_4\big([\underbrace{2 \; 2\; \cdots \; 2}_{j}]_4+n\big)\neq b_4(n)$ when $n=-[\underbrace{1\; 1\; \cdots \; 1}_{j}]_4-1$ and when $n=[\underbrace{1\; 1\; \cdots \; 1}_{j+1}]_4+1$, so this interval cannot be extended.
\end{theorem}

Extensions to the above theorems for higher even bases will be given in Section \ref{sec:higher bases}.

For each of the functions $f_4(n)$ and $b_4(n)$, it is also natural to ask the question: what is the function's rate of growth?  Coons and Spiegelhofer \cite{coons2017maximal} studied a generalized Stern sequence $s_d(n)$, which we will not define here other than to note that it is trivially equal to $f_d(n-1)$.  Defant \cite{defant2016upper}, quoted in Coons and Spiegelhofer, found the maxima for $s_d(n)$ and where they occur. We adapt the statement of this result to show the location of the maxima for $f_d(n)$.

\begin{theorem}[\cite{defant2016upper}, Proposition~2.1]\label{thm:defant}
The maximal value of $f_d(n)$ in the interval $\left[d^{k-2},d^{k-1}\right)$ is $F_k$, the $k$-th Fibonacci number.  Moreover, this maximal value occurs for the first time at the value $n=A_k$, where
\begin{equation*}
A_k=\begin{cases}
\left[(1\; 0)^{\frac{k}{2}-1} \; 0\right]_d,&\mbox{if }k\text{ is even},\\
\left[(1\; 0)^{\frac{k-1}{2}}\right]_d,&\mbox{if }k\text{ is odd}.
\end{cases}
\end{equation*}
\end{theorem}

We show a similar rate of growth for $b_4(n)$.

\begin{theorem}\label{thm:balanced binary maxima}
Let $r>0$.  If $r$ is odd, then $b_4\left(\left[(1\; 2)^{\lceil r/2\rceil}\right]_4\right)=F_{r+3}$.  If $r$ is even, then $b_4\left(\left[(1\; 2)^{r/2} \;2\right]_4\right)=F_{r+3}$.  In both cases, $F_{r+3}$ is the maximal value of $b_4(n)$ in the interval \[ \left[ -[\underbrace{1\; 1\; \cdots \; 1}_{r+1}]_4 , [\underbrace{1\; 1\; \cdots \; 1}_{r+2}]_4 \right] .\]
\end{theorem}

We will provide a partial proof of Theorem \ref{thm:balanced binary maxima} that makes no reference to the $f_4$ function, and only relies on properties of $b_4$. We will then provide a full proof of Theorem \ref{thm:balanced binary maxima} by making use of Theorem \ref{thm:equality on interval} and Theorem \ref{thm:defant}.

This paper is organized as follows.  In Section \ref{sec:recurrence relations}, we prove recurrence relations for $f_4(n)$ and $b_4(n)$ that resemble the recurrences satisfied by the Stern sequence.  The proofs of Theorems \ref{thm:equality on interval}, \ref{thm:no good k}, and \ref{thm:balanced binary maxima} appear in Sections \ref{sec:equality on interval}, \ref{sec:no good k}, and \ref{sec:maxima}, respectively.  Section \ref{sec:higher bases} describes generalizations of Theorems \ref{thm:no good k} and \ref{thm:equality on interval} to general even bases.

\section{Recurrence relations}\label{sec:recurrence relations}

Both hyperquaternary and balanced quaternary representations satisfy Stern-like recurrence relations resembling those in \eqref{eq:stern recurrence}.  The following theorem gives a recurrence relation for hyperquaternary representations.

\begin{theorem}\label{thm:hyperquaternary_recursion}
    For any non-negative integer $n$, we have that
\begin{equation}\label{eq:hyperquaternary_recursion}
        f_4(n)=\begin{cases}
            f_4\left(\frac{n-k}{4}\right),& \text{ if }n\equiv k\pmod{4} \text{ and } k\in\{1,2,3\},\\
            f_4\left(\frac{n}{4}\right)+f_4\left(\frac{n}{4}-1\right),&\text{ if } n\equiv 0 \pmod{4}.
        \end{cases}
    \end{equation}
\end{theorem}

A proof of the above statement via generating functions can be found in \cite{courtright2004arithmetic}. A direct proof via congruence classes, similar to the proof of Theorem \ref{balanced_quaternary_recursion} below, is also possible.

We now give the analogous recurrence relation for balanced quaternary representations.

\begin{theorem}\label{balanced_quaternary_recursion}
For any integer $n$, we have that
\begin{equation*}
    b_4(n)=\begin{cases}
            b_4\left(\frac{n-k}{4}\right),& \text{ if }n\equiv k\pmod{4} \text{ and } k\in\{-1,0,1\},\\
            b_4\left(\frac{n-2}{4}\right)+b_4\left(\frac{n+2}{4}\right),&\text{ if }n\equiv2\pmod{4}.
        \end{cases}
\end{equation*}
\end{theorem}

\begin{proof}
    Let $n$ be an integer, and consider a balanced quaternary representation of $n$ as follows:
    \begin{equation*}
        n=\sum_{i=0}^\infty\delta_i4^i,\quad\delta_i\in\{-2,-1,0,1,2\}.
    \end{equation*}
    Suppose $n\equiv k\pmod{4}$ with $k\in\{-1,0,1\}$.  Then $\delta_0=k$, so we have
    \begin{align*}
        n&=k+\delta_14+\delta_24^2+\delta_34^3+\cdots,\\
        n-k&=4\left(\delta_1+\delta_24+\delta_34^2+\cdots\right),\\
        \frac{n-k}{4}&=\delta_1+\delta_24+\delta_34^2+\cdots,
    \end{align*}
    and so each balanced quaternary representation of $n$ corresponds to a unique balanced quaternary representation of $\frac{n-k}{4}$. Similarly, by reversing the order of the above equations, each balanced quaternary representation of $\frac{n-k}{4}$ corresponds to a unique balanced quaternary representation of $n$. Thus, $b_4(n)=b_4\left(\frac{n-k}{4}\right)$.

    Suppose $n\equiv2\pmod{4}$.  Then $\delta_0=2$ or $\delta_0=-2$.  If $\delta_0=2$, then we have
    \begin{align*}
        n&=2+\delta_14+\delta_24^2+\delta_34^3+\cdots,\\
        n-2&=4\left(\delta_1+\delta_24+\delta_34^2+\cdots\right),\\
        \frac{n-2}{4}&=\delta_1+\delta_24+\delta_34^2+\cdots.
    \end{align*}
    If $\delta_0=-2$, then we have
    \begin{align*}
        n&=-2+\delta_14+\delta_24^2+\delta_34^3+\cdots,\\
        n+2&=4\left(\delta_1+\delta_24+\delta_34^2+\cdots\right),\\
        \frac{n+2}{4}&=\delta_1+\delta_24+\delta_34^2+\cdots.
    \end{align*}
    Therefore, by the same reasoning as above, each balanced quaternary representation of $n$ that has $\delta_0=2$ (respectively, $\delta_0=-2$) corresponds to a unique balanced quaternary representation of $\frac{n-2}{4}$ (respectively, $\frac{n+2}{4}$) and vice versa. Thus,   $b_4(n)=b_4\left(\frac{n-2}{4}\right)+b_4\left(\frac{n+2}{4}\right)$.
\end{proof}

These recurrence relations provide us with a quick proof of the following important fact.

\begin{lemma}\label{lem:function positivity}
    For $n\ge 0$, we have $f_4(n)\ge 1$, and for $n\in \mathbb{Z}$, we have $b_4(n)\ge 1$.
\end{lemma}

In other words, there exists at least one hyperquaternary representation for all non-negative integers and at least one balanced quaternary representation for all integers.

\begin{proof}
    The existence of a  hyperquaternary representation for each non-negative integer follows immediately from the existence of the standard quaternary representation. However, we can also prove this fact by noting that $f_4(0)=f_4(1)=f_4(2)=f_4(3)=1$ and that Theorem \ref{thm:hyperquaternary_recursion} implies that for any $n\in\mathbb{N}$, $n\ge 4$, there exists $m\in\mathbb{N}$ with $0<m<n$ such that $f_4(n)\ge f_4(m)$. We can repeat this latter fact inductively until we arrive at $f_4(n)\ge f_4(m)$, where $m\in \{1,2,3\}$. 

    The proof for $b_4(n)$ proceeds similarly, noting that $b_4(-1)=b_4(0)=b_4(1)=1$ and that Theorem \ref{balanced_quaternary_recursion} shows that for any $n\in\mathbb{Z}$ with $|n|\ge 2$, there exists $m\in\mathbb{Z}$ with $|m|<|n|$ and $b_4(n)\ge b_4(m)$.
    \end{proof}

\section{Proof of Theorem \ref{thm:equality on interval}}\label{sec:equality on interval}

In this section, we prove Theorem \ref{thm:equality on interval}.  To establish the shifted identity \[
f_4\big([\underbrace{2 \; 2\; \cdots \; 2}_{j}]_4+n\big)=b_4(n)
    \] on the given interval, we first identify the largest intervals on which the difference \[
f_4\big([\underbrace{2 \; 2\; \cdots \; 2}_{j}]_4+n\big)-b_4(n)
    \] is equal to zero.  To this end, we let $D_j(n)$ be this difference:
\[
D_j(n) = f_4\big([\underbrace{2 \; 2\; \cdots \; 2}_{j}]_4+n\big)-b_4(n).
\]

Our first goal will be to prove the following theorem, which is the first half of Theorem \ref{thm:equality on interval} rephrased.

\begin{theorem}\label{thm: D_j zero values}
    For $j\ge 1$, we have
    \[
    D_j(n)=0, \qquad n\in I_j:=\left[ -[\underbrace{1\; 1\; \cdots \; 1}_{j}]_4 , [\underbrace{1\; 1\; \cdots \; 1}_{j+1}]_4 \right].
    \]
\end{theorem}

Note that
\[
I'_j:=I_j+[\underbrace{2\; 2\; \cdots \; 2}_j]_4=\left[ [\underbrace{1\; 1\; \cdots \; 1}_{j}]_4 , [1 \underbrace{3\; 3\; \cdots \; 3}_j]_4 \right],
\]
and that $I'_j\cap I'_{j+1}=\left[ [\underbrace{1\; 1\; \cdots \; 1}_{j+1}]_4 , [1 \underbrace{3\; 3\; \cdots \; 3}_j]_4 \right]$ will always be non-empty, and so $\bigcup I'_j=\mathbb{N}$.  In particular, for any $n\in\mathbb{N}$, we can find a $j$ such that $f_4(n)=b_4\big(n-[\underbrace{2\; 2\; \cdots \; 2}_j]_4\big)$.

We begin our analysis of $D_j(n)$ by proving recurrence relations similar to those in Theorems \ref{thm:hyperquaternary_recursion} and \ref{balanced_quaternary_recursion}.
\begin{lemma}\label{lem:Djn recurrence}
    For $j\ge 1$, we have 
    \[
    D_j(n) = \begin{cases}
        D_{j-1}\left(\frac{n-k}{4} \right), & \text{ if }n\equiv k\pmod{4} \text{ and } k\in\{-1,0,1\},\\
        D_{j-1}\left( \frac{n+2}{4}\right) + D_{j-1}\left( \frac{n-2}{4}\right) , & \text{ if } n\equiv 2\pmod{4}
    \end{cases}
    \]
\end{lemma}

\begin{proof}
    First, we suppose that $n\equiv k\pmod{4}$ with $k\in\{-1,0,1\}$. Then 
    $[\underbrace{2 \; 2\; \cdots \; 2}_{j}]_4+n\equiv k'\pmod{4}$ with $k'=k+2$, so $k'\in \{1,2,3\}$. We now apply the definition of $D_j(n)$ together with the recurrence relations for $f_4$ and $b_4$ in Theorems \ref{thm:hyperquaternary_recursion} and \ref{balanced_quaternary_recursion}:
    \begin{align*}
        D_j(n) &=  f_4\big([\underbrace{2 \; 2\; \cdots \; 2}_{j}]_4+n\big)-b_4(n)\\
        &=  f_4\left( \frac{[\underbrace{2 \; 2\; \cdots \; 2}_{j}]_4+n-k'}{4}\right)-b_4\left(\frac{n-k}{4}\right)\\
        &=  f_4\left( \frac{[\underbrace{2 \; 2\; \cdots \; 2}_{j}]_4+n-k-2}{4}\right)-b_4\left(\frac{n-k}{4}\right)\\
        &=  f_4\left([\underbrace{2 \; 2\; \cdots \; 2}_{j-1}]_4+ \frac{n-k}{4}\right)-b_4\left(\frac{n-k}{4}\right)\\
        &= D_{j-1}\left(\frac{n-k}{4}\right),
    \end{align*}
    as desired.

    If $n\equiv 2\pmod{4}$, then $[\underbrace{2 \; 2\; \cdots \; 2}_{j}]_4+n\equiv 0\pmod{4}$, and so again applying the definition of $D_j(n)$ with the recurrence relations for $f_4$ and $b_4$, we obtain the following:
    \begin{align*}
        D_j(n) &=  f_4\big([\underbrace{2 \; 2\; \cdots \; 2}_{j}]_4+n\big)-b_4(n)\\
        &=  f_4\left( \frac{[\underbrace{2 \; 2\; \cdots \; 2}_{j}]_4+n}{4}\right)+f_4\left( \frac{[\underbrace{2 \; 2\; \cdots \; 2}_{j}]_4+n-4}{4}\right)-b_4\left( \frac{n+2}{4}\right)-b_4\left( \frac{n-2}{4}\right)\\
        &= f_4\left( \frac{[\underbrace{2 \; 2\; \cdots \; 2}_{j-1} \; 0]_4+n+2}{4}\right)+f_4\left( \frac{[\underbrace{2 \; 2\; \cdots \; 2}_{j-1} \; 0]_4+n-2}{4}\right)-b_4\left( \frac{n+2}{4}\right)-b_4\left( \frac{n-2}{4}\right)\\
        &= f_4\left( [\underbrace{2 \; 2\; \cdots \; 2}_{j-1}]_4+\frac{n+2}{4}\right)-b_4\left( \frac{n+2}{4}\right)+f_4\left([\underbrace{2 \; 2\; \cdots \; 2}_{j-1}]_4+ \frac{n-2}{4}\right)-b_4\left( \frac{n-2}{4}\right)\\
        &=D_{j-1}\left( \frac{n+2}{4}\right) + D_{j-1}\left( \frac{n-2}{4}\right). \qedhere
    \end{align*}
\end{proof}

We now prove Theorem \ref{thm: D_j zero values} using the recursion for $D_j(n)$.

\begin{proof}[Proof of Theorem \ref{thm: D_j zero values}]

    We will proceed by induction on $j$. For the base case we consider $j=1$. Here we have $D_1(n)=f_4(n+2)-b_4(n)$. Moreover, the interval\[
     \left[ -[\underbrace{1\; 1\; \cdots \; 1}_{j}]_4 , [\underbrace{1\; 1\; \cdots \; 1}_{j+1}]_4 \right]
    \]
    is just $[-1,5]$. We have the following simple values for $f_4$ and $b_4$:

\begin{center}
    \begin{tabular}{c|c|c}
        $n$ & $f_4(n)$ & $b_4(n)$\\
        \hline 
        -2 & & 2\\
        -1 & & 1\\
        0 & 1 & 1\\
        1 &1 & 1 \\
        2 & 1 & 2\\
        3 & 1 & 1\\
        4 & 2 & 1\\
        5 & 1 & 1\\
        6 & 1 & 3\\
        7 & 1\\
        8 & 2\\
        
    \end{tabular}
\end{center}
So we see that $D_1(n)=0$ on $[-1,5]$, and we have the base case. Moreover, $D_1(-2)=-1$ and $D_1(6)=-1$.

Now we will assume that $D_{j-1}(n)=0$ for $n\in I_{j-1}$.

Suppose that $n\equiv k\pmod{4}$ with $k\in\{-1,0,1\}$.  Then $n-k$ is the closest multiple of $4$ to $n$, so if $n\in I_j$, we must have
\[
n-k\in \left[ -[\underbrace{1\; 1\; \cdots \; 1}_{j-1} \; 0]_4 , [\underbrace{1\; 1\; \cdots \; 1}_j\; 0]_4 \right].
\]
Thus, 
\[
\frac{n-k}{4}\in \left[ -[\underbrace{1\; 1\; \cdots \; 1}_{j-1}]_4 , [\underbrace{1\; 1\; \cdots \; 1}_{j}]_4 \right]_4 =I_{j-1}.
\] Thus for these values of $n$, we have, by Lemma \ref{lem:Djn recurrence}, that $D_j(n)=D_{j-1}((n-k)/4)$, and since $(n-k)/4\in I_{j-1}$, our inductive hypothesis shows that $D_{j-1}((n-k)/4)=0$ in this case. 

Now suppose that $n\equiv 2\pmod{4}$. If $n\in I_j$, then in fact, the stronger inclusion,
\[
n\in \left[ -[\underbrace{1\; 1\; \cdots \; 1}_{j}]_4 +3, [\underbrace{1\; 1\; \cdots \; 1}_j\; 0]_4-2 \right]_4 \subset I_j
\]
is true. Therefore, we have that
\begin{align*}
    \frac{n-2}{4} &\in \frac{\left[ -[\underbrace{1\; 1\; \cdots \; 1}_{j}]_4 +3-2, [\underbrace{1\; 1\; \cdots \; 1}_j\; 0]_4-2-2 \right]}{4}\\
    &=  \frac{\left[ -[\underbrace{1\; 1\; \cdots \; 1}_{j}]_4 +1, [\underbrace{1\; 1\; \cdots \; 1}_{j-1}\; 0 \; 0]_4\right]}{4}\\
    &=  \frac{\left[ -[\underbrace{1\; 1\; \cdots \; 1}_{j-1} \; 0]_4 , [\underbrace{1\; 1\; \cdots \; 1}_{j-1}\; 0 \; 0]_4\right]}{4}\\    
    &= \left[ -[\underbrace{1\; 1\; \cdots \; 1}_{j-1} ]_4 , [\underbrace{1\; 1\; \cdots \; 1}_{j-1}\; 0 ]_4 \right] \subset I_{j-1}
\end{align*}
and that
\begin{align*}
    \frac{n+2}{4} &\in \frac{\left[ -[\underbrace{1\; 1\; \cdots \; 1}_{j}]_4 +3+2, [\underbrace{1\; 1\; \cdots \; 1}_{j} \; 0]_4-2+2 \right]}{4}\\
    &=  \frac{\left[ -[\underbrace{1\; 1\; \cdots \; 1}_{j}]_4 +5, [\underbrace{1\; 1\; \cdots \; 1}_{j} \; 0]_4\right]}{4}\\
    &=  \frac{\left[ -[\underbrace{1\; 1\; \cdots \; 1}_{j-2} \; 0\; 0]_4 , [\underbrace{1\; 1\; \cdots \; 1}_{j} \; 0]_4\right]}{4}\\    
    &= \left[ -[\underbrace{1\; 1\; \cdots \; 1}_{j-2}\; 0 ]_4 , [\underbrace{1\; 1\; \cdots \; 1}_{j} ]_4\right] \subset I_{j-1}.
\end{align*}
Thus, for these values of $n$, Lemma \ref{lem:Djn recurrence} tells us that $D_j(n)=D_{j-1}((n-2)/4)+D_{j-1}((n+2)/4)$.  Because we also know that $(n-2)/4,(n+2)/4\in I_{j-1}$, our inductive hypothesis again confirms that $D_j(n)=0$. This completes the proof of Theorem \ref{thm: D_j zero values}.\end{proof}

We emphasize here that in contrast to the binary case, where there exists a $k$ for which $\overline{b}_2(k+n)=f_2(n)$ for all $n\ge 0$, the following lemma shows that our $D_j$ function is not always $0$.

\begin{lemma}\label{lemma:nonzero D_j}
    For any $j\ge 1$, we have that
    \begin{equation}
    D_j([\underbrace{1 \; 1 \; \cdots \;1 }_{j+1}]_4+1)=D_j(-[\underbrace{1 \; 1\; \cdots \; 1}_j]-1)=-1. \label{eq:bad D_j values}
    \end{equation}
    In particular, we cannot extend Theorem \ref{thm: D_j zero values} to be true on any larger interval containing $I_j$.
\end{lemma}

\begin{proof}
    We will prove that $D_j([\underbrace{1 \; 1 \; \cdots \;1 }_{j+1}]_4+1)=-1$ as the proof that $D_j(-[\underbrace{1 \; 1\; \cdots \; 1}_j]-1)=-1$ goes similarly.
    
    Note that as part of our proof of Theorem \ref{thm: D_j zero values}, we showed $D_1([1\; 1]_4+1)=D_1(6)=-1$, which is Equation \eqref{eq:bad D_j values} with $j=1$. We now proceed inductively. We assume that \eqref{eq:bad D_j values} holds for $j-1$ and use this to show it holds for $j$. We have $[\underbrace{1 \; 1 \; \cdots \;1 }_{j+1}]_4+1\equiv 2\pmod{4}$, so Lemma \ref{lem:Djn recurrence} tells us that 
    \begin{align*}
        D_j([\underbrace{1 \; 1 \; \cdots \;1 }_{j+1}]_4+1)&=D_{j-1}\left(\frac{[\underbrace{1 \; 1 \; \cdots \;1 }_{j}\;0]_4+4}{4}\right)+D_{j-1}\left(\frac{[\underbrace{1 \; 1 \; \cdots \;1 }_{j}\;0]_4}{4}\right)\\
        & = D_{j-1}([\underbrace{1 \; 1 \; \cdots \;1 }_{j}]_4+1)+D_{j-1}([\underbrace{1 \; 1 \; \cdots \;1 }_{j}]_4)\\
        & = -1 +0=-1,
    \end{align*}
    where the last equality comes from applying our inductive hypothesis as well as Theorem \ref{thm: D_j zero values}. 
\end{proof}
 
Theorem \ref{thm:equality on interval} follows immediately from Theorem \ref{thm: D_j zero values} and Lemma \ref{lemma:nonzero D_j}.

\section{Proof of Theorem \ref{thm:no good k}}\label{sec:no good k}

In this section, we will show that there is no $k$ such that $f_4(n+k)=b_4(n)$ for all $n\ge -k$. To accomplish this, we will need to know exactly when $f_4(n+k)=1$ and $b_4(n)=1$, which we will do in the following lemma.

\begin{lemma}\label{lemma:when are f and b equal to 1}
    For $n\ge 1$, we have that $f_4(n)=1$ if and only if the standard quaternary representation of $n=[\delta_m \; \delta_{m-1} \; \dots \; \delta_0]_4$, with $\delta_m\neq 0$ contains only the digits $1,2,3$. For $n\in\mathbb{Z}$, we have that $b_4(n)=1$ if and only if there is a balanced quaternary representation $n=[\delta_m \; \delta_{m-1} \; \dots \; \delta_0]_4$ that contains only the digits $-1,0,1$.
\end{lemma}

\begin{proof}
    Suppose $[\delta_m \; \delta_{m-1} \; \dots \; \delta_0]_4>0$. We have by Theorem \ref{thm:hyperquaternary_recursion} that 
    \[
    f_4([\delta_m \; \delta_{m-1} \; \dots \; \delta_0]_4) = f_4([\delta_m \; \delta_{m-1} \; \dots \; \delta_1]_4)
    \]
    if $\delta_0\in\{1,2,3\}$. If $\delta_0=0$, then instead
            \[
    f_4([\delta_m \; \delta_{m-1} \; \dots \; \delta_0]_4) = f_4([\delta_m \; \delta_{m-1} \; \dots \; \delta_1]_4)+f_4([\delta_m \; \delta_{m-1} \; \dots \; \delta_1]_4-1),
    \]
    and since $f_4(n)\ge 1$ for all $n\ge 0$, we have
        \[
    f_4([\delta_m \; \delta_{m-1} \; \dots \; \delta_0]_4) > f_4([\delta_m \; \delta_{m-1} \; \dots \; \delta_1]_4).
    \]
 Iterating the above $m$ times, we see that 
    \[
    f_4([\delta_m \; \delta_{m-1} \; \dots \; \delta_0]_4) = f_4(\delta_m)
    \]
    if and only if $\delta_i\in\{1,2,3\}$ for $i=0,1,\dots, m-1$. Moreover, $f_4(1)=f_4(2)=f_4(3)=1$, which shows now that 
        \[
    f_4([\delta_m \; \delta_{m-1} \; \dots \; \delta_0]_4) = 1
    \] if and only if $\delta_i\in\{1,2,3\}$ for $i=0,1,\dots, m$, as desired. 

    The proof for the $b_4$ case is similar.
\end{proof}

In order to control the digits of both $n$ and $n+k$, we will want to know that there is a way of expanding $k$ in balanced quaternary with specific properties. This will be guaranteed by the following lemma.

\begin{lemma}\label{lemma:nonstandard balanced binary representation}
    For $k\in\mathbb{Z}$, there exists a balanced quaternary representation $k=[\delta_m \; \delta_{m-1} \; \dots \; \delta_0]$ where $\delta_i \neq -2$ for any $i$.
\end{lemma}

\begin{proof}
    We know by Lemma \ref{lem:function positivity} that $b_4(k)\ge 1$, so there must exist at least one representation of $k=[\delta_m\; \delta_{m-1} \; \dots \; \delta_0]_4$ in balanced quaternary. However, there may be some digits that equal $-2$.

    We can perform the following manipulation to a quaternary representation. Since $4a-2=4(a-1)+2$, we have $[a \; -\!2]_4=[(a-1) \; 2]_4$ and more generally
    \[
     [\delta_m \; \delta_{m-1} \; \dots \; \delta_{i+2} \; a \; -\!2 \; \delta_{i-1} \; \dots \; \delta_0]_4= [\delta_m \; \delta_{m-1} \; \dots \; \delta_{i+2} \; (a-1) \; 2 \; \delta_{i-1} \; \dots \; \delta_0]_4.
    \]
    If we suppose that the $-2$ in position $i$ is the left-most $-2$ in the representation, then $a\in \{-1,0,1,2\}$, and so this manipulation will result in another balanced quaternary representation. If $a\in \{0,1,2\}$, then this will completely eliminate the left-most $-2$ from the representation without introducing a new one. However, if $a=-1$, this process replaces a digit $-2$ from position $i$ with a digit $-2$ in position $i+1$. We can now iterate this procedure on the new digit $-2$ at position $i+1$, and repeat as necessary until we are preceded by a digit in $\{0,1,2\}$. So, for $b\in \{0,1,2\}$
    \begin{align*}
    &[\delta_m \; \delta_{m-1} \; \dots \; \delta_{i+2+j} \; b \; \underbrace{-\!1 \; \dots \; -\!1}_j \; -\!2 \; \delta_{i-1} \; \dots \; \delta_0]_4
    \\
    &\qquad = [\delta_m \; \delta_{m-1} \; \dots \; \delta_{i+2+j} \; (b-1) \; \underbrace{2 \; \dots \; 2}_{j+1} \; \delta_{i-1} \; \dots \; \delta_0]_4.
    \end{align*}
    This manipulation now removes the left-most appearance of $-2$ from the representation without introducing another. Note that we may need to prepend a digit of $0$ before the $\delta_m$ if, say, $\delta_m=-2$ itself.

    By repeatedly applying the manipulations in the previous paragraph, we can remove all $-2$'s from a balanced quaternary representation without changing the value of the number, as desired.
\end{proof}

We now prove Theorem \ref{thm:no good k}.

\begin{proof}[Proof of Theorem \ref{thm:no good k}]
    Suppose to the contrary that there does exist $k$ such that $f_4(n+k)=b_4(n)$ for all $n\ge -k$. Let $k$ have a balanced quaternary representation of $[\delta_m \; \delta_{m-1} \; \dots \; \delta_0]_4$, with $\delta_i\in\{-1,0,1,2\}$ for each $i$, which exists by Lemma \ref{lemma:nonstandard balanced binary representation}.

    We claim that $\delta_i=2$ for each $i$. Suppose to the contrary again that there exists some $\delta_i\neq 2$. Let $n\in\mathbb{N}$ be given by
    \[
    n=[1 \; \underbrace{1 \; 1 \; \dots \; 1}_{m-i} \; -\delta_i \; \underbrace{1 \; 1 \; \dots \; 1}_i ]_4.
    \]
    Note that this expresses $n$ as a balanced quaternary representation with all digits belonging to $\{-1,0,1\}$. Therefore $b_4(n)=1$ by Lemma \ref{lemma:when are f and b equal to 1}. Conversely, consider $n+k$. We see that this can be written as 
    \[
    n+k=[1\; (\delta_m\!+\!1) \; (\delta_{m-1}\!+\!1) \; \dots  \; (\delta_{i+1}\!+\!1) \; 0 \; (\delta_{i-1}\!+\!1) \; (\delta_{i-2}\!+\!1) \; \dots \; (\delta_0\!+\!1)]_4.
    \]
    This expresses $n+k$ as a standard quaternary representation since all digits belong to $\{0,1,2,3\}$. Since the leading digit is a $1$ and there is a $0$ at some point after this, we have that $f_4(n+k)>1$ by Lemma \ref{lemma:when are f and b equal to 1}. In particular, $f_4(n+k)>b_4(n)$, contrary to assumption, so $\delta_i=2$.

    Thus $k=[\underbrace{2\; 2\; \dots \; 2}_{m+1}]$, exactly as in the definition of $D_{m+1}(n)$. However, by Lemma \ref{lemma:nonzero D_j}, we know that $D_{m+1}(n)\neq 0$ for some value of $n$. Thus even for this value of $k$, there exists an $n$ such that $f_4(n+k)\neq b_4(n)$. This proves the theorem.
\end{proof}

\section{Maxima for $f_4(n)$ and $b_4(n)$}\label{sec:maxima}

To begin with, we provide a proof of the first part of Theorem \ref{thm:balanced binary maxima}, namely that 
\[
F_{r+3}= \begin{cases}
    b_4\left(\left[(1\; 2)^{\lceil r/2\rceil}\right]_4\right), & \text{if }r\text{ is odd},\\
    b_4\left(\left[(1\; 2)^{r/2} \;2\right]_4\right), & \text{if }r \text{ is even},
\end{cases}
\]
using only facts about the $b_4$ function.

\begin{proof}
Note that if $r>0$ is odd, then
\[
\left[(1\; 2)^{\lceil r/2\rceil}\right]_4=\sum_{k=1}^{\lceil \frac{r}{2}\rceil} 6\cdot 4^{2k-2}=\frac{2}{5}\left(4^{r+1}-1\right),
\]
while if $r>0$ is even, 
\[
\left[(1\; 2)^{r/2} \;2\right]_4=2+\sum_{k=1}^{\frac{r}{2}}6\cdot4^{2k-1}=2+\frac{8}{5}\left(4^r-1\right).
\]
The leftmost expressions in the equalities were used in the statement of Theorem \ref{thm:balanced binary maxima}, but we will use all three types of expressions in the proof and will move between equivalent expressions as needed.  We now proceed with a proof by strong induction.

As base cases, we consider $r=1$ and $r=2$.  When $r=1$, we have $b_4\left(\frac{2}{5}(4^2-1)\right)=b_4(6)=3=F_4$.  When $r=2$, we have $b_4\left(2+\frac{8}{5}\left(4^2-1\right)\right)=b_4(26)=5=F_5$.

Suppose the result holds for every positive integer less than $r$.  We will show the result holds for $r$.

First suppose $r$ is odd.  We must show $b_4\left(\frac{2}{5}\left(4^{r+1}-1\right)\right)=F_{r+3}$.  To apply the recurrence relations, we note 
\begin{align*}
    \frac{2}{5}\left(4^{r+1}-1\right)=\sum_{k=1}^{\lceil r/2\rceil} 6\cdot 4^{2k-2}=6\left(4^0+4^2+4^4+\cdots+4^{2\lceil r/2 \rceil-2}\right)\equiv 2\bmod 4.
\end{align*}
Using the recurrence relation for $b_4(n)$ given in Theorem \ref{balanced_quaternary_recursion}, we have
\begin{align*}
b_4\left(\frac{2}{5}\left(4^{r+1}-1\right)\right)=b_4\left(\frac{1}{4}\left(\frac{2\left(4^{r+1}-1\right)}{5}-2\right)\right)+b_4\left(\frac{1}{4}\left(\frac{2\left(4^{r+1}-1\right)}{5}+2\right)\right)=I + II.
\end{align*}

We begin by examining $I$.  Consider
\begin{align*}
    \frac{\frac{2(4^{r+1}-1)}{5}-2}{4}&=\frac{1}{4}\left(\frac{2(4^{r+1}-1)}{5}-2\right)=\frac{1}{4}\left(\sum_{k=1}^{\lceil\frac{r}{2}\rceil}6\cdot4^{2k-2}-2\right)\\
    &=\frac{1}{4}\left(6(4^0+4^2+4^4+\cdots+4^{2\lceil\frac{r}{2}\rceil-2})-2\right)\\
    &=\frac{1}{4}\left(4+6(4^2+4^4+\cdots+4^{2\lceil\frac{r}{2}\rceil-2})\right)\\
    &=1+6\left(4^1+4^3+\cdots+4^{2\lceil\frac{r}{2}\rceil-3}\right)\\
    &\equiv1\bmod4.
\end{align*}
Thus $I=b_4\left(\frac{\frac{2(4^{r+1}-1)}{5}-2}{4}\right)=b_4\left(\frac{\frac{\frac{2(4^{r+1}-1)}{5}-2}{4}-1}{4}\right)$ by the recursion for $b_4(n)$ given in Theorem \ref{balanced_quaternary_recursion}.

Then
\begin{align*}
I&=b_4\left(\frac{1}{4}\left(\frac{\frac{2(4^{r+1}-1)}{5}-6}{4}\right)\right)=b_4\left(\frac{1}{4}\cdot\frac{1}{4}\left(\frac{2(4^{r+1}-1)}{5}-6\right)\right)\\
&=b_4\left(\frac{1}{4}\cdot\frac{1}{4}\left(\sum_{k=1}^{\lceil\frac{r}{2}\rceil}6\cdot4^{2k-2}-6\right)\right)=b_4\left(\frac{1}{4}\cdot\frac{1}{4}\left(\sum_{k=2}^{\lceil\frac{r}{2}\rceil}6\cdot4^{2k-2}\right)\right)\\
&=b_4\left(\sum_{k=2}^{\lceil\frac{r}{2}\rceil}6\cdot4^{2k-4}\right)
=b_4\left(\sum_{k=1}^{\lceil\frac{r-2}{2}\rceil}6\cdot4^{2k-2}\right)=F_{r+1}
\end{align*}
by the inductive hypothesis.

We know $I=F_{r+1}$, and we will show $II=F_{r+2}$ to conclude $b_4\left(\frac{2}{5}\left(4^{r+1}-1\right)\right)=I+II=F_{r+1}+F_{r+2}=F_{r+3}$.

Since $II=b_4\left(\frac{1}{4}\left(\frac{2\left(4^{r+1}-1\right)}{5}+2\right)\right)$, we begin by considering 
\begin{align*}
   \frac{1}{4}\left(\frac{2\left(4^{r+1}-1\right)}{5}+2\right) &=\frac{1}{4}\left(\sum_{k=1}^{\lceil\frac{r}{2}\rceil}6\cdot4^{2k-2}+2\right)\\
   &=\frac{1}{4}\left(6\left(4^0+4^2+4^4+\cdots+4^{2\lceil\frac{r}{2}\rceil-2}\right)+2\right)\\
   &=\frac{1}{4}\left(8+6\left(4^2+4^4+\cdots+4^{2\lceil\frac{r}{2}\rceil-2}\right)\right)\\
   &=2+6\left(4^1+4^3+\cdots+4^{2\lceil\frac{r}{2}\rceil-3}\right).
\end{align*}
Since $r$ is odd, we know $\lceil\frac{r}{2}\rceil=\lceil\frac{r+1}{2}\rceil$, and $r+1$ is even, so $\lceil\frac{r+1}{2}\rceil=\frac{r+1}{2}$.  Thus
\begin{align*}
    \frac{1}{4}\left(\frac{2\left(4^{r+1}-1\right)}{5}+2\right) &=2+6\left(4^1+4^3+\cdots+4^{2\lceil\frac{r+1}{2}\rceil-3}\right)\\
    &=2+6\left(4^1+4^3+\cdots+4^{2\frac{(r+1)}{2}-3}\right)\\
    &=2+6\left(4^1+4^3+\cdots+4^{2\frac{(r-1)}{2}-1}\right)\\
    &=2+\sum_{k=1}^{\frac{r-1}{2}}6\cdot4^{2k-1},
\end{align*}
and $b_4\left(2+\displaystyle\sum_{k=1}^{\frac{r-1}{2}}6\cdot4^{2k-1}\right)=F_{r+2}$ by the strong inductive hypothesis.

Then $b_4\left(\frac{2}{5}\left(4^{r+1}-1\right)\right)=I+II=F_{r+1}+F_{r+2}=F_{r+3}$, and it remains to consider the case when $r$ is even.

Now suppose $r$ is even.  We must show $b_4\left(2+\frac{8}{5}\left(4^r-1\right)\right)=F_{r+3}$.  To apply the recurrence relation for $b_4(n)$, we note 
\begin{align*}
2+\frac{8}{5}\left(4^r-1\right)&=2+\sum_{k=1}^{\frac{r}{2}}6\cdot 4^{2k-1}\\
&=2+6\left(4^1+4^3+4^5+\cdots 4^{r-1}\right)\\
&\equiv 2\bmod 4.
\end{align*}
Using the recurrence for $b_4(n)$, 
\[
b_4\left(2+\frac{8}{5}\left(4^r-1\right)\right)=b_4\left(\frac{2+\frac{8}{5}(4^r-1)-2}{4}\right)+b_4\left(\frac{2+\frac{8}{5}(4^r-1)+2}{4}\right)=III+IV.
\]
We begin by examining $III$.

First note that since $r$ is even, $\frac{r}{2}=\lceil\frac{r}{2}\rceil=\lceil\frac{r-1}{2}\rceil$.
We have 
\begin{align*}
    \frac{2+\frac{8}{5}(4^r-1)-2}{4}&=\frac{\frac{8}{5}(4^r-1)}{4}=\frac{1}{4}\cdot\frac{8}{5}(4^r-1)\\
    &=\frac{1}{4}\cdot\sum_{k=1}^{\frac{r}{2}}6\cdot4^{2k-1}=\sum_{k=1}^{\lceil\frac{r-1}{2}\rceil}6\cdot 4^{2k-2}.
\end{align*}
Then
\[
III=b_4\left(\frac{2+\frac{8}{5}(4^r-1)-2}{4}\right)=b_4\left(\sum_{k=1}^{\lceil\frac{r-1}{2}\rceil}6\cdot 4^{2k-2}\right)=F_{r+2}
\]
by the strong inductive hypothesis.

Now we examine $IV$.  First note that
\begin{align*}
    \frac{2+\frac{8}{5}(4^r-1)+2}{4}&=\frac{4+\frac{8}{5}(4^r-1)}{4}\\
    &=\frac{1}{4}\left(4+\frac{8}{5}(4^r-1)\right)\\
    &=\frac{1}{4}\left(4+\sum_{k=1}^{\frac{r}{2}}6\cdot4^{2k-1}\right)\\
    &=1+6\left(4^0+4^2+4^4+\cdots+4^{r-2}\right)\\
    &\equiv 3\bmod 4.
\end{align*}
Then
\begin{align*}
    IV&=b_4\left(\frac{2+\frac{8}{5}(4^r-1)+2}{4}\right)=b_4\left(1+\sum_{k=1}^{\frac{r}{2}}6\cdot 4^{2k-2}\right)\\
    &=b_4\left(\frac{1}{4}\left(2+\sum_{k=1}^{\frac{r}{2}}6\cdot 4^{2k-2}\right)\right) \text{by the recurrence for } b_4(n)\\
    &=b_4\left(\frac{1}{4}\left(2+6+\sum_{k=2}^{\frac{r}{2}}6\cdot4^{2k-2}\right)\right)=b_4\left(2+\sum_{k=2}^\frac{r}{2}6\cdot 4^{2k-3}\right)\\
    &=b_4\left(2+\sum_{k=1}^{\frac{r-2}{2}}6\cdot 4^{2k-1}\right)=F_{r+1}
\end{align*}
by the strong inductive hypothesis.

Then $b_4\left(2+\frac{8}{5}(4^r-1)\right)=III+IV=F_{r+2}+F_{r+1}=F_{r+3}$, and this completes the case where $r$ is even.  We have shown that if $r$ is odd, then $b_4\left(\left[(1\; 2)^{\lceil r/2\rceil}\right]_4\right)=F_{r+3}$, and if $r$ is even, then $b_4\left(\left[(1\; 2)^{r/2} \;2\right]_4\right)=F_{r+3}$.  This completes the proof. \end{proof}

We now provide a full proof of Theorem \ref{thm:balanced binary maxima} by appealing to Theorem \ref{thm:equality on interval} and Theorem \ref{thm:defant}.

\begin{proof}
First, consider the case that $r$ is odd.  Since $\left[(1\; 2)^{\lceil r/2\rceil}\right]_4$ is in $I_{r+1}=I_{2\left\lceil r/2\right\rceil}$, we apply Theorem \ref{thm: D_j zero values} with $j=2\left\lceil r/2\right\rceil$ to obtain
\begin{align*}
    b_4\left(\left[(1\; 2)^{\lceil r/2\rceil}\right]_4\right)&=f_4\left(\left[(2\; 2)^{\lceil r/2\rceil}\right]_4+\left[(1\; 2)^{\lceil r/2\rceil}\right]_4\right)\\
    &=f_4\left(\left[(3\; 4)^{\lceil r/2\rceil}\right]_4\right)\\
    &=f_4\left(\left[(1\; 0)^{\lceil r/2\rceil}\; 0\right]_4\right).
\end{align*}
The leading term of $\left[(1\; 0)^{\lceil r/2\rceil}\; 0\right]_4$ is $4^{2\lceil r/2\rceil}$, so we see that $\left[(1\; 0)^{\lceil r/2\rceil}\; 0\right]_4$ is in $\left[4^{2\lceil r/2\rceil},4^{2\lceil r/2\rceil+1}\right)$.  Thus the index $k-2$ from the interval $\left[d^{k-2},d^{k-1}\right)$ in Theorem \ref{thm:defant} satisfies $k-2=2\lceil r/2\rceil$.  Rearranging, we have that $\frac{k}{2}-1=\lceil r/2\rceil$, and it is evident that $k$ is even.  We have that the maximal value of $f_4(n)$ in the interval $\left[4^{k-2},4^{k-1}\right)$ is $f_4\left(\left[(1\; 0)^{\frac{k}{2}-1}\; 0\right]_4\right)=F_k$  from Theorem \ref{thm:defant}; equivalently, the maximal value of $f_4(n)$ in the interval $\left[4^{2\lceil r/2\rceil},4^{2\lceil r/2\rceil+1}\right)$ is $f_4\left(\left[(1\; 0)^{\lceil r/2\rceil}\; 0\right]_4\right)$, which is equal to $F_{2\lceil r/2\rceil+2}=F_{r+3}$. Suppose now that the maximum of $b_4(n)$ in $I_{r+1}$ exceeds $F_{r+3}$ and that this occurs at $m\in I_{r+1}$. Then \[m+[(2\; 2)^{\lceil r/2\rceil}]_4\in I'_{r+1}=\left[ [\underbrace{1\; 1\; \cdots \; 1}_{r+1}]_4 , [1 \underbrace{3\; 3\; \cdots \; 3}_{r+1}]_4 \right].\] But the interval $I'_{r+1}$ is contained in $[4^{r+2}, 4^{r+4})$, and the value of $b_4(m)=f_4(m+[(2\; 2)^{\lceil r/2\rceil}]_4)$ is therefore at most $F_{r+3}$ by Theorem \ref{thm:defant}, which is a contradiction. Therefore, the maximum of $b_4(n)$ in $I_{r+1}$ is $F_{r+3}$ as desired. 

Now, consider the case that $r$ is even.  Since $\left[(1\; 2)^{r/2}\; 2\right]_4$ is in $I_{r+1}$, we apply Theorem \ref{thm: D_j zero values} with $j=r+1$ to obtain
\begin{align*}
    b_4\left(\left[(1\; 2)^{r/2}\; 2\right]_4\right)&=f_4\left(\left[(2\; 2)^{r/2}\; 2\right]_4+\left[(1\; 2)^{r/2}\; 2\right]_4\right)\\
    &=f_4\left(\left[(3\; 4)^{r/2}\; 4\right]_4\right)\\
    &=f_4\left(\left[(1\; 0)^{\frac{r+2}{2}}\right]_4\right).
\end{align*}
The leading term of $\left[(1\; 0)^{\frac{r+2}{2}}\right]_4$ is $4^{r+1}$, so we see that $\left[(1\; 0)^{\frac{r+2}{2}}\right]_4$ is in $\left[4^{r+1},4^{r+2}\right)$.  Thus the index $k-2$ from Theorem \ref{thm:defant} satisfies $k-1=r+1$. 
 Therefore, $k=r+3$ and $\frac{k-1}{2}=\frac{r+2}{2}$.  We have that the maximal value of $f_4(n)$ in the interval $\left[4^{k-2},4^{k-1}\right)$ is $f_4\left(\left[(1\; 0)^{\frac{k-1}{2}}\right]_4\right)=F_k$ from Theorem \ref{thm:defant}; equivalently, the maximal value of $f_4(n)$ in the interval $\left[4^{r+1},4^{r+2}\right)$ is $f_4\left(\left[(1\; 0)^{\frac{r+2}{2}}\right]_4\right)=F_{r+3}$. The proof that  $b_4(n)$ on the interval  $I_{r+1}$ is maximized at $\left[(1\; 2)^{r/2}\; 2\right]_4$ is identical to the odd case.
\end{proof}

\section{Higher even bases}\label{sec:higher bases}

Having restricted our attention in the previous sections to the base-$4$ case, we now consider base $d=2\ell$, where $\ell\in\mathbb{N}$ 
 with $\ell>2$. Our basic definitions extend from base 4 in the expected way. We let $[\delta_m \; \delta_{m-1} \; \dots \; \delta_0]_{2\ell}$ denote $\sum \delta_i (2\ell)^i$, let $f_{2\ell}(n)$ denote the number of ways of writing $n$ in hyper-$2\ell$-ary (with digits $\delta_i\in \{0,1,2,\dots, 2\ell\}$), and let $b_{2\ell}(n)$ denote the number of ways of writing $n$ in balanced $2\ell$-ary (with digits $\delta_i\in \{-\ell, -\ell+1,\dots, \ell-1,\ell\}$).

The following recurrence relations follow by similar proofs to the base-4 case:
    for any positive integer $n\equiv k\bmod 2\ell$ with $0\le k < 2\ell$, we have that
    \begin{equation}\label{eq: hyper dary recursion}
        f_{2\ell}(n)=\begin{cases}
            f_{2\ell}\left(\frac{n-k}{2\ell}\right),&\mbox{if }k\neq0,\\
            f_{2\ell}\left(\frac{n}{2\ell}\right)+f_{2\ell}\left(\frac{n}{2\ell}-1\right),&\mbox{if }k=0.
        \end{cases}
    \end{equation}
Similarly,     for any integer $n\equiv k\bmod 2\ell$ with $-\ell< k \le \ell$, we have that
    \begin{equation}\label{eq: balanced dary recursion}
        b_{2\ell}(n)=\begin{cases}
            b_{2\ell}\left(\frac{n-k}{2\ell}\right),&\mbox{if }k \neq \ell,\\
            b_{2\ell}\left(\frac{n-\ell}{2\ell}\right)+b_{2\ell}\left(\frac{n+\ell}{2\ell}\right),&\mbox{if }k=\ell.
        \end{cases}
    \end{equation}
A proof of \eqref{eq: hyper dary recursion} can again be found in \cite{courtright2004arithmetic}, but proofs of both \eqref{eq: hyper dary recursion} and \eqref{eq: balanced dary recursion} also follow by a recursive method, as used in our proof of Theorem \ref{balanced_quaternary_recursion}.

We have the following generalization of Theorems \ref{thm:no good k} and \ref{thm:equality on interval}.

\begin{theorem}\label{thm:main theorem general case}
    There does not exist an integer $k$ such that $f_{2\ell}(n+k)=b_{2\ell}(n)$ for all $n\ge -k$. However, for each $j\in \mathbb{N}$, we have that
    \[
     f_{2\ell}\big([\underbrace{\ell \; \ell\; \cdots \; \ell}_{j}]_{2\ell}+n\big)=b_{2\ell}(n)\] for \begin{equation} n \in \left[ -[\underbrace{(\ell-1)\; (\ell-1)\; \cdots \; (\ell-1)}_{j}]_{2\ell} , [\underbrace{(\ell-1)\; (\ell-1)\; \cdots \; (\ell-1)}_{j+1}]_{2\ell} \right]. \label{eq:general n interval}
    \end{equation}
    Moreover, $f_{2\ell}\big([\underbrace{\ell \; \ell\; \cdots \; \ell}_{j}]_{2\ell}+n\big)\neq b_{2\ell}(n)$ when $n$ equals $-[\underbrace{(\ell-1)\; (\ell-1)\; \cdots \; (\ell-1)}_{j}]_{2\ell}-1$ or $[\underbrace{(\ell-1)\; (\ell-1)\; \cdots \; (\ell-1)}_{j+1}]_{2\ell}+1$, so this interval cannot be extended.
\end{theorem}

The proof of this more general statement follows by a nearly exact proof to the base-4 case. As such, we shall elide over most of it.

We let $D_{j,2\ell}(n) = f_{2\ell}\big([\underbrace{\ell \; \ell \; \cdots \; \ell}_{j}]_{2\ell}+n\big)-b_{2\ell}(n)$. The next statement, which is an analogue of Lemma \ref{lem:Djn recurrence}, follows by an identical proof to the base-$4$ case:
\[
D_{j,2\ell}(n) = \begin{cases}
    D_{j-1,2\ell}\left( \frac{n-k}{2\ell}\right), &  \text{ if }n\not\equiv \ell\pmod{2\ell},\\
        D_{j-1,2\ell}\left( \frac{n+\ell}{2\ell}\right) + D_{j-1,2\ell}\left( \frac{n-\ell}{2\ell}\right) , & \text{ if } n\equiv \ell\pmod{2\ell}.
\end{cases}
\]
The analogue of Theorem \ref{thm: D_j zero values} is that $D_{j,2\ell}(n)=0$ for $n$ satisfying \eqref{eq:general n interval}. The proof will again be an inductive argument. In the base-$4$ case, we proved the base case by simple examination, but we will need a more refined argument for the base case in the general base-$2\ell$ case. To be clear, the base case requires us to show that \[
D_{1,2\ell}(n) = f_{2\ell}\left( \ell+n\right) - b_{2\ell}(n)=0
\]
for \[
n\in[-[\ell-1]_{2\ell},[(\ell-1)\; (\ell-1)]_{2\ell}].
\]

\iffalse
Let $D_j(n) = f_{2\ell}\big([\underbrace{\ell \; \ell \; \cdots \; \ell}_{j}]_{2\ell}+n\big)-b_{2\ell}(n)$. 
\textcolor{red}{Do we want to make this $D_{j,2\ell}(n)$ or is that too much?} To generalize Theorem \ref{thm: D_j zero values}, we need to show that \textcolor{red}{some words here} holds for the base case $D_1(n)$. Namely, we want to show that
\[
D_1(n) = f_{2\ell}\left( \ell+n\right) - b_{2\ell}(n)=0
\]
for \[
n\in[-[\ell-1]_{2\ell},[(\ell-1)\; (\ell-1)]_{2\ell}].
\]
\fi

First, we note that $f_{2\ell}(0)=b_{2\ell}(0)=1$, as the only representation of $0$ in either hyper-$2\ell$-ary or balanced $2\ell$-ary is $[0]_{2\ell}$. Using \eqref{eq: hyper dary recursion} and \eqref{eq: balanced dary recursion}, we see that $f_{2\ell}(n)=1$ for $n\in \{0,1,\dots, 2\ell-1\}$ and $b_{2\ell}(n)=1$ for $n\in \{-\ell+1,-\ell+2, \dots, \ell-1\}$ and also that $f_{2\ell}(2\ell)=f_{2\ell}(0)+f_{2\ell}(1)=2$ and $b_{2\ell}(\ell) = b_{2\ell}(0)+b_{2\ell}(1)=2$. Thus far we have shown that $D_{1,2\ell}(n)=0$ for  $n\in [-[\ell-1]_{2\ell},[\ell]_{2\ell}]$.

Now suppose that $n$ is of the form $2\ell m+k$ for $m=0,1,\dots,\ell-1$ and $k\in [-(\ell-1),\ell-1]$. Then, by the first case of \eqref{eq: balanced dary recursion}, we have that $b_{2\ell}(n) = b_{2\ell}(m)=1$. Similarly, by the first case of \eqref{eq: hyper dary recursion}, $f_{2\ell}(\ell+n)=f_{2\ell}(m)=1$. 

On the other hand, suppose that $n$ is of the form $2\ell m +\ell$ with $0\le m \le \ell-2$. Then $b_{2\ell}(n) = b_{2\ell}(m)+b_{2\ell}(m+1) = 2$. Similarly, $f_{2\ell}(n+\ell) = f_{2\ell}(m)+f_{2\ell}(m+1)=2$. Combined, the previous two paragraphs show that $D_{1,2\ell}(n)=0$ for $n\in[-[\ell-1]_{2\ell},[(\ell-1)\; (\ell-1)]_{2\ell}].$

Before proceeding to the inductive step, we wish to show one more aspect of the case where $j=1$.  If $n=-\ell$, then $b_{2\ell}(n) = b_{2\ell}(0)+b_{2\ell}(-1) = 2$, but $f_{2\ell}(n+\ell) = f_{2\ell}(0)=1$, so $D_{1,2\ell}(-\ell) = -1$. If $n=2\ell(\ell-1)+\ell$, then we have that
\begin{align*}
b_{2\ell}(2\ell(\ell-1)+\ell) &= b_{2\ell}(\ell-1)+b_{2\ell}(\ell) \\
&=b_{2\ell}(\ell-1)+b_{2\ell}(0)+b_{2\ell}(1)\\
&= 3
\end{align*}
and $f_{2\ell}(n+\ell)=f_{2\ell}(2\ell(\ell))=2$.  Then $D_{1,2\ell}(2\ell(\ell-1)+\ell)=-1$. This shows the $j=1$ case of the $2\ell$-analogue of Lemma \ref{lemma:nonzero D_j}. 

The inductive steps of Theorem \ref{thm: D_j zero values} and Lemma \ref{lemma:nonzero D_j} follow by analogous methods.

It remains to show that there does not exist an integer $k$ such that $f_{2\ell}(n+k)=b_{2\ell}(n)$ for all $n\ge -k$. All the lemmas and proofs of Section \ref{sec:no good k} extend naturally. For completion, we mention how a few details are changed. 

First, recall that in Lemma \ref{lemma:nonstandard balanced binary representation}, we wanted to show that for any $k\in\mathbb{Z}$, there is a balanced quaternary representation of $k$ that does not use the digit $-2$.  We use made of the relation $[a\; -2]_4 = [(a-1) \; 2]_4$. Now we want to show that for any $k\in\mathbb{Z}$, there is a balanced $2\ell$-ary representation of $k$ that does not use the digit $\ell$.  Here we can make use of the easily verifiable relation $[ a\; -\ell]_{2\ell}= [(a-1)\; \ell]_{2\ell}$ to eliminate all copies of $-\ell$ from a representation. 

Second, as part of the proof of Theorem \ref{thm:no good k}, we let $k=[\delta_m \; \delta_{m-1}\; \dots \;\delta_0]_4$ with $\delta_i \in \{-1,0,1,2\}$ for each $i$. If $\delta_i=2$ for some specific $i$, then we chose
\[
    n=[1 \; \underbrace{1 \; 1 \; \dots \; 1}_{m-i} \; -\delta_i \; \underbrace{1 \; 1 \; \dots \; 1}_i ]_4.
\]
In the general base-$2\ell$ case, we now want to let $k=[\delta_m \; \delta_{m-1}\; \dots \;\delta_0]_{2\ell}$ with $\delta_i \in \{-(\ell-1),-(\ell-2),\dots , \ell-1,\ell\}$ for each $i$. If $\delta_i=\ell$ for some specific $i$, then we choose
\[
n=[1 \; \underbrace{(\ell-1) \; (\ell-1) \; \dots \; (\ell-1)}_{m-i} \; -\delta_i \; \underbrace{(\ell-1) \; (\ell-1) \; \dots \; (\ell-1)}_i ]_{2\ell}.
\]
This guarantees that the standard $2\ell$-ary representation  of $n+k$ has a leading $1$ and contains at least one $0$, which is necessary for the proof of Theorem \ref{thm:main theorem general case}.

This completes the catalog of changes necessary to prove Theorem \ref{thm:main theorem general case}.

\end{document}